\documentclass[12pt,a4paper]{amsart}
\usepackage{amssymb}
\usepackage{mathrsfs}
\headsep=1cm \numberwithin{equation}{section}
\newtheorem{thm}{Theorem}[section]
\newtheorem{cor}[thm]{Corollary}
\newtheorem{lem}[thm]{Lemma}

\newtheorem{prop}[thm]{Proposition}
\theoremstyle{definition}
\newtheorem{defn}[thm]{Definition}
\theoremstyle{remark}
\newtheorem{rem}[thm]{Remark}
\numberwithin{equation}{section}

\newcommand\Supp{\operatorname{Supp}}
\newcommand\Ass{\operatorname{Ass}}
\newcommand\mAss{\operatorname{mAss}}

\newcommand\Ann{\operatorname{Ann}}

\newcommand\Rad{\operatorname{Rad}}

\newcommand\Ext{\operatorname{Ext}}

\newcommand\Att{\operatorname{Att}}

\begin{document}\title[Annihilators of local cohomology]
{Cohomological dimension filtration and annihilators of top local cohomology modules}
\author{Ali Atazadeh, Monireh Sedghi$^*$ and Reza Naghipour
\vspace*{0.2cm}}
\address{Department of Mathematics, Azarbaijan Shahid Madani University, Tabriz, Iran.}
\email{aalzp2002@yahoo.com}

\address{Department of Mathematics, Azarbaijan Shahid Madani University, Tabriz, Iran.}
\email {m\_sedghi@tabrizu.ac.ir}
\email {sedghi@azaruniv.ac.ir}

\address{Department of Mathematics, University of Tabriz, Tabriz, Iran;
and School of Mathematics, Institute for Research in Fundamental
Sciences (IPM), P.O. Box 19395-5746, Tehran, Iran.}
\email{naghipour@ipm.ir} \email {naghipour@tabrizu.ac.ir}
\thanks{ 
This research  was in part supported by a grant from IPM.\\
$^*$Corresponding author: e-mail: {\it m\_sedghi@tabrizu.ac.ir} (Monireh Sedghi)}%

\begin{abstract}
Let $\frak a$ denote an ideal in a commutative Noetherian ring $R$ and $M$ a finitely generated $R$-module. In this paper, we introduce the concept of the cohomological dimension filtration $\mathscr{M} =\{M_i\}_{i=0}^c$,  where $ c={\rm cd} (\frak a,M)$ and $M_i$ denotes the largest submodule of $M$ such that ${\rm cd} (\frak a,M_i)\leq i.$ Some properties of this filtration are investigated. In particular, in the case that $(R, \frak m)$ is local and $c= \dim M$, we are able to determine the annihilator of the top local cohomology module $H_{\frak a}^c(M)$. In fact, it is shown that ${\rm Ann}_R(H_{\frak a}^c(M))= {\rm Ann}_R(M/M_{c-1}).$ As a consequence, it follows that there exists an ideal $\frak b$ of $R$ such that $\Ann_R(H_{\frak a}^{c}(M))=\Ann_R(M/H_{\frak b}^{0}(M))$. This generalizes the main results of \cite{BAG} and \cite{Lyn}.
\end{abstract}
\maketitle
\section{Introduction}
Let $R$ be an arbitrary commutative Noetherian ring
(with identity), $\frak a$ an ideal of $R$ and let $M$ be a finitely generated $R$-module. An important problem concerning local cohomology is determining the annihilators of the $i^{\rm th}$ local cohomology module $H_{\frak a}^i(M)$. This problem has been studied by several authors; see for example \cite{HK}, \cite{Lyn}, \cite{Ly1}, \cite{NA},  \cite{NC},  \cite{Sc1},   and  has led to some interesting results. In particular,
in \cite{BAG} Bahmanpour et al., proved an interesting result about the annihilator of the $d$-th local cohomology module $\Ann_R(H_{\frak m}^d(M))$,  in the case $(R, \frak m)$ is a complete local ring.

The purpose of the present paper is to introduce the concept of the cohomological dimension filtration $\mathscr{M} =\{M_i\}_{i=0}^c$,  where
$c={\rm cd} (\frak a,M)$ and $M_i$ denotes the largest submodule of $M$ such that ${\rm cd} (\frak a,M_i)\leq i.$  Because $M$  is a Noetherian $R$-module,
it follows easily from \cite[Theorem 2.2]{DNT} that the submodules $M_i$ are well-defined. They also form increasing family of submodules. Some properties of this filtration are investigated. In particular, in the case that $(R, \frak m)$ is local and $c= \dim M$, we are able to determine the annihilator of the top local cohomology module $H_{\frak a}^c(M)$. In fact, it is shown that $$\Ann_R(H_{\frak a}^c(M))= \Ann_R(M/M_{c-1}).$$ As a consequence, it follows that there exists an ideal $\frak b$ of $R$ such that $$\Ann_R(H_{\frak a}^{c}(M))=\Ann_R(M/H_{\frak b}^{0}(M)).$$ This generalizes the main results of \cite{Ba}, \cite{BAG} and \cite{Lyn}.

As a main result in the second section, we describe in more details the structure of the cohomological dimension filtration $\mathscr{M} =\{M_i\}_{i=0}^c$,
in terms of the reduced primary decomposition of $0$ in $M$. Namely, if  $0=\cap_{j=1}^nN_j$ denotes a reduced primary decomposition of the zero submodule $0$ in $M$ such that $N_j$ is a $\frak p_j$-primary submodule of $M$ and ${\frak a}_i:=\Pi_{{\rm cd}(\frak a, R/{\frak p_j})\leq i}\frak p_j$,
we shall show that:

\begin{thm}
Let $R$ be a Noetherian ring and $\frak a$ an ideal of $R$. Let $M$ be a non-zero finitely generated $R$-module with finite cohomological dimension $c:={\rm cd}(\frak a, M)$ with respect to $\frak a$ and let  $\mathscr{M}=\{ M_i \}_{i=0}^c$ be the cohomological dimension filtration of $M$. Then, for all integers $0\leq i \leq c$,

${\rm (i)}$ $M_i=H_{{\frak a}_i}^0(M)=\cap_{{\rm cd}(\frak a, R/{\frak p_j})>i}N_j;$

$\rm(ii)$ ${\rm Ass}_RM_i=\{\frak p\in {\rm Ass}_RM|\,{\rm cd}(\frak a, R/\frak p)\leq i\};$

$\rm(iii)$ ${\rm Ass}_RM/M_i=\{\frak p\in {\rm Ass}_RM|\,{\rm cd}(\frak a, R/\frak p)> i\};$

$\rm(iv)$ ${\rm Ass}_RM_i/M_{i-1}=\{\frak p\in {\rm Ass}_RM|\,{\rm cd}(\frak a, R/\frak p)= i\}.$
\end{thm}

Pursuing this point of view further we establish some
results about the annihilator of the top local cohomology modules.  More precisely, as a main result of the third section,
we derive the following consequence of Theorem {\rm1.1},  which will be describe the  annihilator of the top local cohomology module
$H_{\frak a}^{\dim M}(M)$.

\begin{thm}
Let $\frak a$ denote an ideal of a local (Noetherian) ring $R$ and let $M$ be a  finitely
generated $R$-module of dimension $c$ such that $H^c_{\frak a}(M)\neq 0$. Then
$$\Ann_R(H_{\frak a}^{c}(M))=\Ann_R(M/M_{c-1}).$$
\end{thm}
Several corollaries of this result are given.  A typical result in this direction is the following, which is a generalization of the main results of \cite[Theorem 2.6]{BAG}, \cite[Theorem 3.2]{Ba} and \cite[Theorem 2.4]{Lyn} for an  ideal $\frak a$ in an arbitrary local  ring $R$.
\begin{cor}
Let $R$ be a local (Noetherian) ring and $\frak a$ an ideal of $R$. Let $M$ be a non-zero finitely generated $R$-module of dimension $d$ such that $H_{\frak a}^d(M)\neq 0$. Then $$\Ann_R(H_{\frak a}^{d}(M))=\Ann_R(M/H_{\frak b}^{0}(M))=\Ann_R(M/\cap_{{\rm cd}(\frak a, R/{\frak p_j})=d}N_j).$$ Here $0=\cap_{j=1}^nN_j$ denotes a reduced primary decomposition of the zero submodule $0$ in $M$ and $N_j$ is a $\frak p_j$-primary submodule of $M$, for all $j=1,\dots ,n$ and $\frak b:=\Pi_{{\rm cd}(\frak a, R/{\frak p_j})\neq d}\frak p_j.$
\end{cor}

Throughout this paper, $R$ will always be a commutative Noetherian
ring with non-zero identity and $\frak a$ will be an ideal of $R$. For any $R$-module $L$, the
$i^{\rm th}$ local cohomology module of $L$ with support in $V(\frak a)$
is defined by
$$H^i_{\frak a}(L) := \underset{n\geq1} {\varinjlim}\,\, \Ext^i_R(R/\frak a^n, L).$$

For each $R$-module $L$, we denote by
 ${\rm Assh}_RL$ (resp. $\mAss_RL$) the set $\{\frak p\in \Ass
_RL:\, \dim R/\frak p= \dim L\}$ (resp. the set of minimal primes of
$\Ass_RL$).  For an Artinian $R$-module $A$,  we shall use $\Att_R A$ to denote the set of attached prime
ideals of $A$.
Also, for any ideal $\frak a$ of $R$, we denote $\{\frak p \in {\rm Spec}\,R:\, \frak p\supseteq \frak a \}$ by
$V(\frak a)$. In addition, for any ideal $\frak{b}$ of $R$, {\it the radical} of $\frak{b}$, denoted by $\Rad(\frak{b})$, is defined to
be the set $\{x\in R \,: \, x^n \in \frak{b}$ for some $n \in \mathbb{N}\}$.  Finally, if $(R, \mathfrak{m})$ is a local (Noetherian)  ring and
$M$ a finitely generated $R$-module, then $\hat{R}$ (resp. $\hat{M}$)
denotes the completion of $R$ (resp.  $M$) with respect to the
$\mathfrak{m}$-adic topology. For any unexplained notation and terminology we refer
the reader to \cite{BS} and \cite{Mat}.\\

\section{Cohomological dimension filtration}

For an $R$-module $M$, the {\it cohomological dimension} of $M$ with respect to an ideal $\frak a$ of $R$ is defined as $${\rm cd}(\frak a, M):= \sup\{i\in \mathbb{Z}|\, H^i_{\frak a}(M)\neq 0\}.$$ If $(R, \frak m)$ is local and $\frak a=\frak m$, it is known that ${\rm cd}(\frak a, M)= \dim M.$

The purpose of this section is to introduce the notion of {\it cohomological dimension filtration} (abbreviated as cd-filtration) of $M$, which is a generalization of the concept of {\it dimension filtration} introduced by P. Schenzel \cite{Sc2}. Specifically, let $\frak a$ be an ideal of $R$ and $M$ a non-zero finitely generated $R$-module with finite cohomological dimension $c:=cd(\frak a, M)$.  For an integer $0\leq i \leq c$, let $M_i$ denote the largest submodule of $M$ such that ${\rm cd}(\frak a, M_i) \leq i$. Because of the maximal condition of a Noetherian $R$-module, it follows easily from \cite[Theorem 2.2]{DNT} that the submodules $M_i$ of $M$ are well-defined. Moreover, it is clear that $M_{i-1} \subseteq M_i$ for all $1\leq i \leq c$.
\begin{defn}
The increasing filtration $\mathscr{M}=\{ M_i \}_{i=0}^c$ of submodules of $M$ is called the
{\it  cohomological dimension filtration} (abbreviated as cd-filtration) of $M$, where $ c={\rm cd} (\frak a,M)$.
\end{defn}
Before investigating some properties of the cohomological dimension filtration, we state the following lemma which plays a key role in this paper.

\begin{lem}
Let $R$ be a Noetherian ring and $\frak a$ an ideal of $R$. Let $M$ and $N$ be two finitely generated
$R$-modules such that ${\rm Supp}N \subseteq {\rm Supp}M$. Then $${\rm cd}(\frak a, N) \leq {\rm cd}(\frak a, M).$$

\end{lem}

\proof See \cite[Theorem 2.2]{DNT}.\qed\\

\begin{prop}
Let $R$ be a Noetherian ring and $\frak a$ an ideal of $R$. Let $M$ be a non-zero finitely generated $R$-module with finite cohomological dimension $c:={\rm cd}(\frak a, M)$  and let  $\mathscr{M}=\{ M_i \}_{i=0}^c$ be the cd-filtration of $M$. Then, for all integers $0\leq i \leq c$, we have $$M_i=H_{{\frak a}_i}^0(M)=\cap_{{\rm cd}(\frak a, R/{\frak p_j})>i}N_j.$$ Here $0=\cap_{j=1}^nN_j$ denotes a reduced primary decomposition of the zero submodule $0$ in $M$, $N_j$ is a $\frak p_j$-primary submodule of $M$ and $a_i=\Pi_{{\rm cd}(\frak a, R/{\frak p_j})\leq i}\frak p_j.$
\end{prop}
\proof First, we show the equality $H_{{\frak a}_i}^0(M)=\cap_{{\rm cd}(\frak a, R/{\frak p_j})>i}N_j.$  To do this, the inclusion $\cap_{{\rm cd}(\frak a, R/{\frak p_j})>i}N_j \subseteq H_{{\frak a}_i}^0(M)$ follows by easy arguments about the primary decomposition of the zero submodule $0$ of $M$. In order to prove the opposite inclusion, suppose, the contrary is true. Then, there exists $x\in H_{{\frak a}_i}^0(M)$ such that $x\not\in \cap_{{\rm cd}(\frak a, R/{\frak p_j})>i}N_j.$ Thus there exists an integer $t$ such that $x\not\in N_t$ and ${\rm cd}(\frak a, R/{\frak p_t})>i.$ Now, as $x\in H_{{\frak a}_i}^0(M)$, it follows that there is an integer $s_i \geq 1$ such that $\frak a_i^{s_i}x=0$, and so $\frak a_i^{s_i}x \subseteq N_t$. Because of $x\not\in N_t$ and $N_t$ is a $\frak p_t$-primary submodule, it yields that $\frak a_i \subseteq \frak p_t$. Hence there is an integer $j$ such that $\frak p_j \subseteq \frak p_t$ and ${\rm cd}(\frak a, R/{\frak p_j}) \leq i.$ Therefore, in view of Lemma {\rm2.2}, we have $${\rm cd}(\frak a, R/{\frak p_t}) \leq {\rm cd}(\frak a, R/{\frak p_j}) \leq i,$$ which is a contradiction. Now we show that $M_i=H_{{\frak a}_i}^0(M).$ To do this, let $x\in M_i$. Then, in view of Lemma {\rm2.2}, ${\rm cd}(\frak a, Rx) \leq i$. Now, let $\frak p$ be an arbitrary minimal prime ideal over ${\rm Ann}_RRx$. Then, using Lemma {\rm2.2}, we see that ${\rm cd}(\frak a, R/\frak p) \leq i$. On the other hand, since $\frak p \in {\rm Ass}_RRx$, it follows that $\frak p \in {\rm Ass}_RM$, and so there is $1\leq j \leq n$ such that $\frak p_j=\frak p$.
 Hence we have
$${\frak a}_i\subseteq \cap_{{\rm cd}(\frak a, R/{\frak p_j})\leq i}\frak p_j\subseteq \cap_{\frak p\in \mAss_R(Rx)}\frak p=\Rad({\rm Ann}_RRx).$$
Therefore, there exists an integer $n_i\geq1$ such that $\frak a_i^{n_i}\subseteq {\rm Ann}_RRx$, and hence $\frak a_i^{n_i}x=0$. That is $x\in H_{{\frak a}_i}^0(M)$, and so $M_i\subseteq H_{{\frak a}_i}^0(M).$ On the other hand, as ${\rm Supp}\,H_{\frak a_i}^0(M)\subseteq V(\frak a_i),$ it follows that for every $\frak p\in {\rm Supp}\,H_{\frak a_i}^0(M)$, there exists an integer $j\geq 1$ such that $\frak p_j \subseteq \frak p$ and ${\rm cd}(\frak a, R/{\frak p_j})\leq i.$ Since ${\rm Supp}\,R/\frak p\subseteq {\rm Supp}\,R/\frak p_j$, it follows from Lemma {\rm2.2} that ${\rm cd}(\frak a, R/{\frak p})\leq {\rm cd}(\frak a, R/{\frak p_j})\leq i.$ Therefore, in view of \cite[Corollary 2.2]{Sc3}, we have ${\rm cd}(\frak a, H_{\frak a_i}^0(M))\leq i$. Now, it follows from the maximality of $M_i$  that $M_i=H_{{\frak a}_i}^0(M)$, as required.\qed\\
\begin{defn}
Let $R$ be a Noetherian ring and $\frak a$ an ideal of $R$. Let $M$ be a non-zero finitely generated $R$-module with finite cohomological dimension $c:={\rm cd}(\frak a, M)$. We denote by $T_R(\frak a, M)$ the largest submodule of $M$ such that ${\rm cd}(\frak a, T_R(\frak a, M))< c$.
\end{defn}

Using Lemma {\rm2.2}, it is easy to check that $T_R(\frak a, M)=\cup \{ N\leqslant M|\, {\rm cd}(\frak a, N) < c\}.$ In particular, for a local ring $(R, \frak m)$, we denote $T_R(\frak m, M)$ by $T_R(M)$. Thus $$T_R(M)=\cup \{ N\leqslant M|\, {\rm dim}N<{\rm dim}M\}.$$
\begin{rem}
Let $R$ be a Noetherian ring, $\frak a$ an ideal of $R$ and $M$ a non-zero finitely generated $R$-module with finite cohomological dimension $c:={\rm cd}(\frak a, M)$. Let $\{M_i \}_{i=0}^c$ be a cd-filtration of $M$. Then $T_R(\frak a,M)=M_{c-1}$ and by Proposition {\rm2.3}, we have $$T_R(\frak a,M)=H_{\frak b}^0(M)=\cap_{{\rm cd}(\frak a, R/{\frak p_j})=c}N_j,$$ where $0=\cap_{j=1}^nN_j$ denotes a reduced primary decomposition of the zero submodule $0$ in $M$, $N_j$ is a $\frak p_j$-primary submodule of $M$ and $\frak b=\Pi_{{\rm cd}(\frak a, R/{\frak p_j})\neq c}\frak p_j.$
\end{rem}
The next proposition provides an information about the associated primes of the cohomological dimension filtration of $M$.
\begin{prop}
Let $R$ be a Noetherian ring and $\frak a$ an ideal of $R$. Let $M$ be a non-zero finitely generated $R$-module with finite cohomological dimension $c:={\rm cd}(\frak a, M)$, and let $\{M_i \}_{i=0}^c$ be the cd-filtration of $M$. Then, for all integers $0\leq i \leq c$,

$\rm(i)$ ${\rm Ass}_RM_i=\{\frak p\in {\rm Ass}_RM|\,{\rm cd}(\frak a, R/\frak p)\leq i\},$

$\rm(ii)$ ${\rm Ass}_RM/M_i=\{\frak p\in {\rm Ass}_RM|\,{\rm cd}(\frak a, R/\frak p)> i\},$ and

$\rm(iii)$ ${\rm Ass}_RM_i/M_{i-1}=\{\frak p\in {\rm Ass}_RM|\,{\rm cd}(\frak a, R/\frak p)= i\}.$

\end{prop}
\proof In view of  Proposition {\rm2.3}, $M_i=H_{{\frak a}_i}^0(M)$, and so by \cite[Section 2.1, Proposition 10]{Bo}, we have
$$\Ass_RM_i=\Ass_R\,M\cap V(\frak a_i).$$
Now, $\rm(i)$ follows from  Lemma {\rm2.2}. In order to show $\rm(ii)$, use \cite[Exercise 2.1.12]{BS}. Finally, for proving $\rm(iii)$, as $M_i/M_{i-1}\subseteq M/M_{i-1}$, it follows  from part (ii) that ${\rm Ass}_RM_i/M_{i-1}\subseteq {\rm Ass}_RM$ and ${\rm cd}(\frak a, R/\frak p)\geq i$, for every $\frak p\in {\rm Ass}_RM_i/M_{i-1}$. Furthermore, in view of the exact sequence
 $$0 \longrightarrow M_{i-1} \longrightarrow M_i \longrightarrow M_i/M_{i-1} \longrightarrow 0,$$
 and  Lemma {\rm2.2}, we have $${\rm cd}(\frak a, M_i/M_{i-1}) \leq {\rm cd}(\frak a, M_i)\leq i.$$ Hence, using again Lemma {\rm2.2}, we obtain that ${\rm cd}(\frak a, R/\frak p)\leq i$, for all $\frak p\in {\rm Ass}_RM_i/M_{i-1}.$ Therefore, $${\rm Ass}_RM_i/M_{i-1}\subseteq \{\frak p\in {\rm Ass}_RM|\,{\rm cd}(\frak a, R/\frak p)= i\}.$$ Now, let $\frak p \in {\rm Ass}_RM$ and ${\rm cd}(\frak a, R/\frak p)= i.$ Then, by virtue of part $\rm(i)$, $\frak p \in {\rm Ass}_RM_i.$ As $\frak p \notin {\rm Ass}_RM_{i-1}$, it follows from the exact sequence, $$0 \longrightarrow M_{i-1} \longrightarrow M_i \longrightarrow M_i/M_{i-1} \longrightarrow 0$$ that $\frak p \in {\rm Ass}_RM_i/M_{i-1}$, and so $${\rm Ass}_RM_i/M_{i-1}=\{\frak p\in {\rm Ass}_RM|\,{\rm cd}(\frak a, R/\frak p)= i\}.$$ \qed\\

\section{Annihilators of top local cohomology modules}

The main point of this section is to determine the annihilator of top local cohomology modules in terms of the reduced primary decomposition of
the zero submodule. Our main result is Theorem 3.5. The following lemmas and proposition play a key role in the proof of that theorem.

\begin{lem}{\rm(cf.\cite[Lemma 7.3.1]{BS})}
 Let $R$ be a Noetherian ring and $\frak a$ an ideal of $R$. Let $M$ be a non-zero finitely generated $R$-module with finite dimension $d$ such that  $H_{\frak a}^d(M)\neq 0.$ Set $G:=M/T_R(\frak a, M).$ Then

$\rm(i)$ ${\rm cd}(\frak a,G)=d,$

$\rm(ii)$ $G$ has no non-zero submodule of cohomological dimension (with respect to $\frak a$) less than $d$,

$\rm(iii)$ ${\rm Ass}_RG={\rm Att}_RH_{\frak a}^d(G)=\{\frak p\in {\rm Ass}_RM|\,{\rm cd}(\frak a, R/\frak p)= d \},$

$\rm(iv)$ $H_{\frak a}^d(G)\cong H_{\frak a}^d(M).$

\end{lem}
\proof The assertion follows easily from  Proposition {\rm2.6}$\rm(ii)$,  \cite[Theorem 2.5]{Di}, Lemma 2.2 and  the exact sequence
$$0 \longrightarrow T_R(\frak a, M)\longrightarrow M \longrightarrow G\longrightarrow 0.$$\qed\\

Before bringing the next lemma let us recall the important notion of a cofinite  module with respect to an ideal. For an ideal $\frak a$ of $R$, an $R$-module
$M$ is said to be $\frak a$-cofinite if $M$ has support in $V (\frak a)$ and $\Ext^i_R(R/\frak a, M)$ is finitely generated for each $i$. The
concept of $\frak a$-cofinite  modules was introduced by R. Hartshorne \cite{Ha}.

\begin{lem}
Let $(R, \frak m)$ be a local (Noetherian) ring such that $\hat{R}$ is integral over $R$. Let $\frak a$ be an ideal of $R$ and $M$ a non-zero
finitely generated $R$-module of dimension $d$. Then
$$\Att_R(H_{\frak a}^{d}(M))=\{\frak p\in {\rm Assh}_RM|\, \Rad(\frak a+\frak p)=\frak m\}.$$
\end{lem}
\proof Let $\frak p\in \Att_R(H_{\frak a}^{d}(M))$. Then it follows from \cite[Theorem A]{DY} that $\frak p\in {\rm Assh}_RM$
and ${\rm cd}(\frak a, R/\frak p)=d$. Now, since by \cite[Theorem 3]{DM} the $R$-module $H_{\frak a}^{d}(M)$ is $\frak a$-cofinite,
it yields from  \cite[Theorem 2.2]{AB} that $\Rad(\frak a+\frak p)=\frak m$.  Hence
$$\Att_R(H_{\frak a}^{d}(M))\subseteq\{\frak p\in {\rm Assh}_RM|\, \Rad(\frak a+\frak p)=\frak m\}.$$
To prove the reverse inclusion, let $\frak p\in {\rm Assh}_RM$ be such that $\Rad(\frak a+\frak p)=\frak m$. Then, as $\dim R/\frak p=d$ and
$$H^d_{\frak a}(R/\frak p)\cong H^d_{\frak a(R/\frak p)}(R/\frak p)\cong H^d_{\frak m}(R/\frak p) ,$$
it follows that ${\rm cd}(\frak a, R/\frak p)=d$, and so in view of  \cite[Theorem A]{DY}, $\frak p\in \Att_R(H_{\frak a}^{d}(M))$.
This completes the proof. \qed\\

\begin{cor}
Let $(R, \frak m)$ be a complete local (Noetherian) ring,  $\frak a$ an ideal of $R$ and $M$ a non-zero
finitely generated $R$-module of dimension $d$. Then
$$\Att_R(H_{\frak a}^{d}(M))=\{\frak p\in {\rm Assh}_RM|\, \Rad(\frak a+\frak p)=\frak m\}.$$
\end{cor}
\proof The result follows from Lemma 3.2. \qed\\

The following proposition  will serve to shorten the proof of the main theorem.
\begin{prop}
Let $(R, \frak m)$ be a complete local (Noetherian) ring and $\frak a$ an ideal of $R$. Let $M$ be a non-zero finitely generated $R$-module of dimension $d$ such that $H_{\frak a}^d(M)\neq 0$. Then $$\Ann_R(H_{\frak a}^{d}(M))=\Ann_R(M/T_R(\frak a,M)).$$

\end{prop}
\proof Let $G:=M/T_R(\frak a,M).$ Then, in view of Lemma 3.1, it is enough for us to show that $\Ann_R(H_{\frak a}^{d}(G))=\Ann_R(G).$ To do this,  it follows easily from Lemma 3.1(iii) and Corollary 3.3  that  $\frak m=\Rad(\frak a+\Ann_R(G)).$ Consequently, $H_{\frak a}^{d}(G)\cong H_{\frak m}^{d}(G)$, and hence, in view of \cite[Theorem 2.6]{BAG}, we have $$\Ann_R(H_{\frak a}^{d}(G))=\Ann_R(G/T_R(G)).$$ Now, since $${\rm cd}(\frak a,T_R(G))\leq {\rm dim}\,T_R(G)<{\rm dim}\,G,$$ it follows from Lemma 3.1(ii) that $T_R(G)=0$, and so $\Ann_R(H_{\frak a}^{d}(G))=\Ann_R(G),$ as required. \qed\\

We are now ready to  prove the main theorem of this
section, which generalizes all of the previous results concerning
the annihilators of the top local cohomology modules.
\begin{thm}
Let $(R, \frak m)$ be a local (Noetherian) ring and $\frak a$ an ideal of $R$. Let $M$ be a non-zero finitely generated $R$-module of dimension $d$ such that $H_{\frak a}^d(M)\neq 0$. Then $$\Ann_R(H_{\frak a}^{d}(M))=\Ann_R(M/T_R(\frak a,M)).$$
\end{thm}
\proof In view of Lemma 3.1, we may assume that $T_R(\frak a,M)=0.$ Now, as $$\Ann_R(M)\subseteq \Ann_R(H_{\frak a}^{d}(M)),$$ it is enough to show that $\Ann_R(H_{\frak a}^{d}(M))\subseteq \Ann_R(M).$ To this end, let $x\in \Ann_R(H_{\frak a}^{d}(M))$ and we show that $xM=0$. Suppose the contrary, that $xM\neq 0$. Then, as $T_R(\frak a, M)=0$, it follows that ${\rm cd}(\frak a, xM)=d.$ Hence ${\rm cd}({\frak a}\hat{R}, x\hat{M})=d,$ and so
$xH_{\frak a \hat{R}}^d(\hat{M})\neq 0$. Because, if $xH_{\frak a \hat{R}}^d(\hat{M})=0$,  then
 $x\hat{R}\subseteq \Ann_{\hat{R}}(H_{\frak a\hat{R}}^d(\hat{M})).$  Hence, in view of Proposition 3.4, $$x\hat{R}\subseteq \Ann_{\hat{R}}(\hat{M}/T_{\hat{R}}(\frak a\hat{R}, \hat{M})),$$ and so $x\hat{M}\subseteq T_{\hat{R}}(\frak a\hat{R}, \hat{M}).$ Therefore,
 ${\rm cd}(\frak a\hat{R}, x\hat{M})<d,$ which is a contradiction. Consequently, $xH_{\frak a}^d(M)\neq 0$. That is
 $x\not\in \Ann_R(H_{\frak a}^{d}(M))$, which is a contradiction. \qed\\

 The first application of Theorem 3.5 improves a result of Coung et al. in  \cite[Lemma 3.2]{CDN}.
 \begin{cor}
Let $R$ be a local (Noetherian) ring and $\frak a$ an ideal of $R$. Let $M$ be a non-zero finitely generated $R$-module of dimension $d$ such that $H_{\frak a}^d(M)\neq 0$.Then $V(\Ann_R(H_{\frak a}^d(M)))=\Supp (M/T_R(\frak a, M))$.
\end{cor}
\proof In view of Theorem 3.5,  we have $$V(\Ann_R(H_{\frak a}^d(M))=V(\Ann_R(M/T_R(\frak a, M)))=\Supp (M/T_R(\frak a, M)),$$  as required.\qed\\

\begin{cor}
Let $R$ be a local (Noetherian) ring and $\frak a$ an ideal of $R$ such that $H_{\frak a}^{{\rm dim}R}(R)\neq 0$. Then $\Ann_R(H_{\frak a}^{{\rm dim}R}(R))$ is the largest ideal of $R$ such that $${\rm cd}\,(\frak a, \Ann_R(H_{\frak a}^{{\rm dim}R}(R))) < {\rm dim}R.$$
\end{cor}
\proof The assertion follows from Theorem 3.5. \qed\\
\begin{prop}
Let $R$ be a local (Noetherian) ring and $\frak a$ an ideal of $R$. Let $M$ be a non-zero finitely generated $R$-module of dimension $d$ such that $H_{\frak a}^d(M)\neq 0$. Then $$\Ann_R(H_{\frak a}^{d}(M))=\Ann_R(M/H_{\frak b}^{0}(M))=\Ann_R(M/\cap_{{\rm cd}(\frak a, R/{\frak p_j})=d}N_j).$$ Here $0=\cap_{j=1}^nN_j$ denotes a reduced primary decomposition of the zero submodule $0$ in $M$, $N_j$ is a $\frak p_j$-primary submodule of $M$, for all $j=1,\dots ,n$ and $\frak b=\Pi_{{\rm cd}(\frak a, R/{\frak p_j})\neq d}\frak p_j.$
\end{prop}
\proof The assertion follows easily from Theorem 3.5 and Remark 2.5.\qed\\

The following corollary is a generalization of the main result of \cite[Theorem 2.4]{Lyn}.

\begin{cor}
Let $R$ be a local (Noetherian) ring of dimension $d$ and $\frak a$ an ideal of $R$ such that $H_{\frak a}^d(R)\neq 0$. Then
 $$\Ann_R(H_{\frak a}^{d}(R))=H_{\frak b}^{0}(R)=\cap_{{\rm cd}(\frak a, R/{\frak p_j})=d}\frak q_j,$$ where $0=\cap_{j=1}^n\frak q_j$ is a reduced primary decomposition of the zero ideal of $R$, $\frak q_j$ is a $\frak p_j$-primary ideal of $R$, for all $1\leq j\leq n$ and
 $\frak b=\Pi_{{\rm cd}(\frak a, R/{\frak p_j})\neq d}\frak p_j.$
\end{cor}
\proof The result follows readily from Corollary 3.8. \qed\\

\begin{prop}
Let $R$ be a local (Noetherian) ring and $\frak a$ an ideal of $R$. Let $M$ be a non-zero finitely generated $R$-module of dimension $d$ such that $H_{\frak a}^d(M)\neq 0$. Then

$\rm(i)$ $\Rad(\Ann_R(H_{\frak a}^d(M))= \bigcap_{\frak p\in \Ass_RM, \, {\rm cd}(\frak a, R/\frak p)=d}\frak p.$

$\rm(ii)$ $\Supp(H_{\frak a}^d(M))\subseteq \bigcup_{\frak p\in \Ass_RM, \, {\rm cd}(\frak a, R/\frak p)=d} V(\frak p+\frak a)$.
\end{prop}
\proof The assertion (i) follows from Proposition 3.8. In order to prove (ii), by using part (i) we have
\begin{eqnarray*}
\Supp(H_{\frak a}^d(M))\subseteq V(\Ann_R(H_{\frak a}^d(M))&=& V(\Rad(\Ann_R(H_{\frak a}^d(M)) \\
&=&V(\cap_{\frak p\in \Ass_RM, \, {\rm cd}(\frak a, R/\frak p)=d}\frak p)\\&=&\cup_{\frak p\in \Ass_RM, \, {\rm cd}(\frak a, R/\frak p)=d} V(\frak p).
\end{eqnarray*}
Now, as $\Supp(H_{\frak a}^d(M))\subseteq V(\frak a)$, it follows that
$$\Supp(H_{\frak a}^d(M))\subseteq (\cup_{\frak p\in \Ass_RM, \, {\rm cd}(\frak a, R/\frak p)=d} V(\frak p))\cap V(\frak a),$$
and thus the desired result follows from this. \qed\\

\begin{cor}
Let $R$ be a local (Noetherian) ring,  $\frak a$ an ideal of $R$  and $x\in R$. Let $M$ be a non-zero finitely generated $R$-module of dimension $d$ such that $H_{\frak a}^d(M)\neq 0$. Then $H_{\frak a}^d(xM)=0$ if and only if $xH_{\frak a}^d(M)=0$. In particular, $\Ann_R(H_{\frak a}^d(M))=0$ if and only if ${\rm cd}(\frak a, rM)=d$, for every non-zero element $r$ of $R$.
\end{cor}
\proof  The assertion follows  readily from Theorem 3.5.\qed\\

The following result is a generalization of \cite[Corollary 2.5]{Lyn} and \cite[Corollary 2.9]{BAG}.\\

\begin{cor}
Let $(R,\frak m)$ be a local (Noetherian) ring of dimension $d$ and $\frak a$ an ideal of $R$.Then the following conditions are equivalent:

$\rm(i)$ $\Ann_RH_{\frak a}^d(R)=0,$

$\rm(ii)$ $\Ass_RR=\Att_RH_{\frak a}^d(R).$
\end{cor}
\proof  In order to show that $\rm(i)\Longrightarrow\rm(ii)$, let $\Ann_R(H_{\frak a}^d(R))=0$ and that $\frak p\in \Ass_RR$. Then $R/{\frak p}\cong Rx$ for some $x(\neq0)\in R$. Thus it follows from Corollary 3.11 that  ${\rm cd}(\frak a, R/\frak p)={\rm cd}(\frak a, Rx)=d$, and so
by \cite[Theorem A]{DY}, $\frak p\in \Att_RH_{\frak a}^d(R)$,  as required.

To prove the implication $\rm(ii)\Longrightarrow\rm(i)$,  in view of Theorem 3.5 and Corollary 3.11, it is enough to show that ${\rm cd}(\frak a, Rx)=d$ for every non-zero element  $x$ of $R$. To do this, in view of \cite[Corollary 2.2]{Sc3},  there exists $\frak p\in \Ass_RRx$ such that
${\rm cd}(\frak a, R/\frak p)={\rm cd}(\frak a, Rx)$.  By virtue of assumption (ii), $\frak p\in \Att_RH_{\frak a}^d(R)$, and so ${\rm cd}(\frak a, R/\frak p)=d$.
Therefore ${\rm cd}(\frak a, Rx)=d$, as required.  \qed\\

\begin{cor}
Let $(R,\frak m)$ be a local (Noetherian) domain of dimension $d$ and $\frak a$ an ideal of $R$ such that $H_{\frak a}^d(R)\neq 0$. Then $\Ann_R(H_{\frak a}^d(R))=0$.
\end{cor}
\proof Since $\Ass_RR={0}$, the assertion follows immediately from Corollary 3.12.\qed\\
\begin{cor}
Let $R$ be a Noetherian domain and $\frak a$ an ideal of $R$ with $ {\rm ht}\, \frak a=n$. Then $\Ann_R(H_{\frak a}^n(R))=0$.
\end{cor}
\proof Suppose the contrary, that $\Ann_R(H_{\frak a}^n(R))\neq 0$. Then there exists a non-zero element $r$ in $\Ann_R(H_{\frak a}^n(R))$.
Hence, we have $rH_{\frak a}^n(R)=0$.  Now, let $\frak q$ be a minimal prime ideal of $\frak a$ such that $ {\rm ht}\, \frak q=n$. Then
$R_{\frak q}$ is a local (Noetherian) domain of dimension $n$ and $rH_{\frak qR_{\frak q}}^n(R_{\frak q})=0$.
Thus $r/1(\neq0)\in \Ann_{R_{\frak q}}(H_{\frak qR_{\frak q}}^n(R_{\frak q}))$, and so by Corollary 3.13, we achive a contradiction.\qed
\begin{cor}
Let $(R, \frak m)$ be a local (Noetherian) ring of dimension $d$ and $\frak a$ an ideal of $R$ such that ${\rm grade}\, \frak a=d$. Then $\Ass_R R=\Att_RH_{\frak a}^d(R).$
\end{cor}
\proof The assertion follows from \cite[Theorem 3.3]{Lyn} and Corollary 3.12.\qed\\
\begin{cor}
Let $R$ be a local (Noetherian) ring and $\frak a$ an ideal of $R$. Let $M$ be a non-zero finitely generated $R$-module of dimension $d$ such that $\Ass_RM=\Att_RH_{\frak a}^d(M)$. Then $\Ann_R(H_{\frak a}^d(M))=\Ann_RM$.
\end{cor}
\proof  Let $0=\cap_{j=1}^nN_j$ denote a reduced primary decomposition of zero submodule $0$ in $M$ such that
$N_j$ is a $\frak p_j$-primary submodule of $M$, for all $j=1,\dots ,n$. Then,
as $\Ass_RM=\Att_RH_{\frak a}^d(M)$ it follows that $\cap_{{\rm cd}(\frak a, R/{\frak p_j})=d}N_j=0$, and so  by Proposition 3.8
$\Ann_R(H_{\frak a}^d(M))=\Ann_RM$, as required. \qed\\

\begin{center}
{\bf Acknowledgments}
\end{center}
The authors would like to thank Professor Peter Schenzel for helpful suggestions, and also the Institute for Research in Fundamental Sciences (IPM), for the financial support.

\end{document}